\documentclass[11pt, a4paper]{article}
\usepackage{enumerate}
\usepackage{amsmath,amsfonts,verbatim,afterpage,theorem,euscript,mathrsfs,amssymb}
\usepackage{amsfonts}
\usepackage{amssymb}
\usepackage{array}
\usepackage{dsfont}
\newtheorem{Definition}{Definition}[section]
\newtheorem{Proposition}{Proposition}[section]
\newtheorem{Lemma}{Lemma}[section]
\newtheorem{Theorem}{Theorem}

\newtheorem{Remark}{Remark}[section]
\usepackage{hyperref}

\def\bkP{{\rm I\kern-.17em P}}
\def\bkC{{\rm I\kern-.5em C}} 

\begin{document}

\title{Real Interpolation method, Lorentz spaces and refined Sobolev inequalities}
\author{Diego CHAMORRO\footnote{$^\dagger$ Laboratoire d'Analyse et de Probabilit\'es, Universit\'e d'Evry Val d'Essonne, 23 Boulevard de France, 91037 Evry Cedex - France} and Pierre-Gilles LEMARI\'E--RIEUSSET$^{\dagger}$}
\maketitle
\abstract{In this article we give a straightforward proof of refined inequalities between Lorentz spaces and Besov spaces and we generalize previous results of H. Bahouri and A. Cohen \cite{COHB}. Our approach is based in the characterization of Lorentz spaces as real interpolation spaces. We will also study the sharpness and optimality of these inequalities.\\[5mm]
}
\section{Introduction}

This paper is a generalization of recent results by H. Bahouri and A. Cohen \cite{COHB} on Lorentz spaces and refined Sobolev inequalities. Let us recall the setting: if $q>1$, the classical Sobolev inequalities $\|f\|_{L^p}\leq C \|f\|_{\dot{W}^{s,q}}$ with $1/p=1/q-s/n$, had been refined by P. G\'erard, Y. Meyer and F. Oru \cite{GERM} using a Besov space in the right-hand side of the inequality:
\begin{equation}\label{Ineq1}
\|f\|_{L^p} \leq C \|f\|_{\dot{W}^{s,q}}^{q/p}\|f\|_{\dot{B}^{s-n/q,\infty}_{\infty}}^{1-q/p}.
\end{equation}
Similarly the inequality $\|f\|_{L^p}\leq C \|f\|_{\dot{B}^{s,q}_q}$ for $q\geq 1$ may be refined by 
\begin{equation}\label{Ineq2}
\|f\|_{L^p} \leq C \|f\|_{\dot{B}^{s,q}_q}^{q/p}\|f\|_{\dot{B}^{s-n/q,\infty}_{\infty}}^{1-q/p}.
\end{equation}
In the previous formula, and for all the following theorems and inequalities, since $\alpha, \beta>0$, we will say that $f\in \dot{B}^{\alpha,p_0}_{q_0}\cap \dot{B}^{-\beta,p_1}_{q_1}$ if $f$ can be writen using the Littlewood-Paley decomposition $f=\displaystyle{\sum_{j\in \mathbb{Z}}}\Delta_j f$ and if the semi-norms $\|\cdot\|_{\dot{B}^{\alpha,p_0}_{q_0}}$ and $\|\cdot\|_{\dot{B}^{-\beta,p_1}_{q_1}}$ are bounded: this way we will choose a natural representation in such spaces.\\[5mm]

In\cite{COHB}, H. Bahouri and A. Cohen show that it is possible to improve the estimate 
$$\|f\|_{L^{p,q}}\leq C \|f\|_{\dot{B}^{s,q}_q}$$
where $L^{p,q}$ is a Lorentz space, into the following inequality
\begin{equation}\label{Ineq3}
\hspace{2cm} \|f\|_{L^{p,q}}\leq C \|f\|_{\dot{B}^{s,q}_q}^{q/p}\|f\|_{\dot{B}^{s-n/q,\infty}_q}^{1-q/p}
\end{equation}
and they prove that this estimate is sharp since it is not possible to replace the norm $\dot{B}^{s-n/q,\infty}_q$ above by a weaker Besov norm $\dot{B}^{s-n/q,\infty}_r$ with $r>q$. They also ask the following question: it is possible to improve the inequality 
$$\|f\|_{L^{p,r}}\leq C \|f\|_{\dot{B}^{s,q}_q}$$
into
\begin{equation}\label{Ineq4}
\hspace{4.4cm} \|f\|_{L^{p,r}}\leq C \|f\|_{\dot{B}^{s,q}_r}^{q/p}\|f\|_{\dot{B}^{s-n/q,\infty}_r}^{1-q/p} \qquad \mbox{when $r\neq q$?}
\end{equation}
The proof they gave for (\ref{Ineq3}) involves difficult estimations of Lorentz norms and could not easily be extended to the case $r\neq q$.\\

In this article we are going to provide an elementary proof of (\ref{Ineq3}), show that the inequality (\ref{Ineq4}) is valid and we will also study a wider family of related inequalities. The main idea is, in the spirit of \cite{LEMA}, to make a systematic use of the characterization of Lorentz spaces as interpolation spaces and to never use the traditional definition of Lorentz spaces as Orlicz spaces.\\

Our first theorem deals with the two first inequalities given above in a quite general form. As we shall see in the proof given in section \ref{Proofs}, this is just a variant of Hedberg's inequality \cite{HEDB}. 
\begin{Theorem}\label{Theorem1}
Let $\alpha, \beta >0$ and $q_0,q_1\in [1+\infty]$. Let $\theta={\alpha\over \alpha+ \beta}\in (0,1)$ and let ${1\over p}={1-\theta\over q_0}+{\theta\over q_1}$. Then there exists a constant $C_0$ such that, for every $f\in  \dot F^{\alpha,q_0}_\infty \cap \dot F^{-\beta,q_1}_\infty(\mathbb{R}^n)$, we have
\begin{equation}\label{Ineq5}
|f(x)| \leq C_0  \bigg(\sup_{j\in\mathbb{Z}} 2^{j\alpha} |\Delta_jf(x)| \bigg)^{1-\theta} \bigg(\sup_{j\in\mathbb{Z}} 2^{-j\beta} |\Delta_jf(x)| \bigg)^{\theta}
\end{equation}
In particular, we get :
\begin{equation}\label{Ineq6}
\|f\|_{L^p}\leq C_0 \|f\|_{\dot F^{\alpha,q_0}_\infty}^{1-\theta}\|f\|_{\dot F^{-\beta,q_1}_\infty}^\theta
\end{equation}
\end{Theorem}
Let us show briefly how to obtain the two first inequalities from this theorem. We observe that, for $0<s<n/q$ and for $1<q<+\infty$, we have $\dot{W}^{s,q}=\dot{F}^{s,q}_2\subset \dot{F}^{s,q}_{\infty}\subset \dot{F}^{s-\frac{n}{q},\infty}_\infty$. Thus, if we set $\alpha=s>0$, $-\beta=s-n/q<0$, $q_0=q$ and $q_1=+\infty$, we deduce from (\ref{Ineq6}) the inequality (\ref{Ineq1}). In the same way, for $0<s<n/q$ and $1\leq q\leq +\infty$, we have $\dot{B}^{s,q}_q=\dot{F}^{s,q}_q\subset \dot{F}^{s,q}_{\infty}\subset \dot{F}^{s-\frac{n}{q},\infty}_\infty$ and we get inequality (\ref{Ineq2}) from inequality (\ref{Ineq6}).\\

Note also that we obtain by similar arguments the following useful inequality:
\begin{equation}\label{Ineq7}
\|f\|_{L^p}\leq C_0 \|f\|_{\dot B^{\alpha,q_0}_{q_0}}^{1-\theta}\|f\|_{\dot B^{-\beta,q_1}_{q_1}}^\theta
\end{equation}
with  $\frac{1}{p}=\frac{1-\theta}{q_0}+\frac{\theta}{q_1}$ and $\theta=\frac{\alpha}{\alpha+\beta}\in ]0,1[$.\\

Our next result studies inequalities (\ref{Ineq3}) and (\ref{Ineq4}). The proof of this theorem reduces to a few lines once we have in mind the characterization of Lorentz spaces as interpolation spaces.
\begin{Theorem}\label{Theorem2}
Let $\alpha, \beta>0$, $q_0,q_1\in [1+\infty]$ with $q_0\neq q_1$. Let $\theta={\alpha\over \alpha+\beta}\in ]0,1[$, ${1\over p}={1-\theta\over q_0}+{\theta\over q_1}$ and let  $r \in [1,+\infty]$.  If $f\in \dot{B}^{\alpha, q_0}_r\cap  \dot{B}^{-\beta, q_1}_r(\mathbb{R}^n)$ then $f\in L^{p,r}(\mathbb{R}^n)$ and we have 
\begin{equation}\label{Ineq8}
\|f\|_{L^{p,r}}\leq C_0 \|f\|_{\dot B^{\alpha,q_0}_r}^{1-\theta}\|f\|_{\dot B^{-\beta,q_1}_r}^\theta.
\end{equation}
\end{Theorem}
In the expression above we have the same index $r$ in Lorentz and Besov spaces, so the the next step is to deal with general Besov spaces $ \dot B^{\alpha,q_0}_{r_0}$ and $\dot B^{-\beta,q_1}_{r_1}$ and to try to control a Lorentz norm in $L^{p,r}$. For this we define $\theta={\alpha\over \alpha+\beta}\in ]0,1[$ and we investigate the validity of the inequality
\begin{equation}\label{Ineq9}
\|f\|_{L^{p,r}}\leq C_0 \|f\|_{\dot B^{\alpha,q_0}_{r_0}}^{1-\theta}\|f\|_{\dot B^{-\beta,q_1}_{r_1}}^\theta.
\end{equation}
A scaling argument\footnote{The norms $\|\ \|_X$ involved in (\ref{Ineq9}) are homogeneous: $\lambda>0\mapsto \|f(\lambda x)\|_X$ is an homogeneous function of $\lambda$.} gives us that necessarily we have  ${1\over p}={1-\theta\over q_0}+{\theta\over q_1}$. For the index $r$, we will have a similar condition as it is explained in the next theorem.\\

It is important to observe that inequality (\ref{Ineq9}) can be studied from the point of view of interpolation theory by the following equivalent problem:
$$\big[\dot B^{\alpha,q_0}_{r_0}, \dot B^{-\beta,q_1}_{r_1}\big]_{\theta, 1}\subset L^{p,r}.$$
However the interpolation between these spaces with these parameters is a delicate issue as it is explained in \cite{KREP} and, to the best of our knowledge, it was not treated before.\\

The next result studies the validity of inequality (\ref{Ineq9}) in some particular cases.

\begin{Theorem}\label{Theorem3}
Let $\alpha,\beta,r^*>0$ and $\theta={\alpha\over \alpha+\beta}$ with ${1\over p}={1-\theta\over q_0}+{\theta\over q_1}$ and ${1\over r^*}={1-\theta\over r_0}+{\theta\over r_1}$. 
\begin{itemize}
\item[1)] For $f\in \dot B^{\alpha,q_0}_{r_0}\cap \dot B^{-\beta,q_1}_{r_1}(\mathbb{R}^n)$ and if $r>r^*$, we have $f\in L^{p,r}(\mathbb{R}^n)$ and
$$\|f\|_{L^{p,r}}\leq C \|f\|_{\dot B^{\alpha,q_0}_{r_0}}^{1-\theta}\|f\|_{\dot B^{-\beta,q_1}_{r_1}}^\theta.$$
\item[2)] Moreover, this inequality is valid for $r=r^*$ in the following cases :
\begin{itemize}
\item[a)] $r=r_0=r_1$,
\item[b)] $r_0=q_0$ and $r_1=q_1$,
\item[c)] $1<p\leq 2$ and $r^*=p$.
\end{itemize}
\item[3)] Finally, the condition $r\geq r^*$ is sharp.
\end{itemize}
\end{Theorem}
As we shall see in section \ref{Proofs}, theorems \ref{Theorem2} \& \ref{Theorem3} are obtained by direct interpolation. However, it is possible to go one step further with the following results:
\begin{Theorem}\label{Theorem4}
Let $\alpha, \beta>0$, $q_0,q_1\in [1+\infty]$, $q_0< q_1$. Let $\theta={\alpha\over \alpha+\beta}\in ]0,1[$ and let ${1\over p}={1-\theta\over q_0}+{\theta\over q_1}$. Let  $q_0\leq r_0\leq r_1\leq q_1$ and let ${1\over r}= {1-\theta\over r_0}+{\theta\over r_1}$. Then we have 
$$ \|f\|_{L^{p,r}}\leq C_0 \|f\|_{\dot B^{\alpha,q_0}_{r_0}}^{1-\theta}\|f\|_{\dot B^{-\beta,q_1}_{r_1}}^\theta.$$
\end{Theorem}
As we shall see, theorem \ref{Theorem4} can be obtained by the use of inequalities (\ref{Ineq7}) and (\ref{Ineq8}) following the ideas given in theorem \ref{Theorem2}.\\

Note now that in all these inequalities we have $\frac{1}{p}=\frac{1-\theta}{q_0}+\frac{\theta}{q_1}$ for $1\leq q_0<q_1\leq +\infty$ and a similar condition for $r, r_0$ and $r_1$: namely $\frac{1}{r}=\frac{1-\theta}{r_0}+\frac{\theta}{r_1}$. We also assumed the following relationship between these parameters:
$$q_0\leq r_0\leq r_1\leq q_1.$$ 
A more general result is given with the next theorem:
\begin{Theorem}\label{Theorem5}
Let $\alpha, \beta>0$, $q_0,q_1\in [1+\infty]$, $q_0\neq q_1$. Let $\theta={\alpha\over \alpha+\beta}\in ]0,1[$ and let ${1\over p}={1-\theta\over q_0}+{\theta\over q_1}$. If $p\leq 2$, we define
\begin{equation*}
I=\bigg\{ (x,y)\in [0,1]^2\   /\ x(1-\theta)+\theta y={1\over p} \bigg\}=[(x_0,y_0),(x_1,y_1)]
\end{equation*}
 with $x_0<y_0$ and $x_1<y_1$, then if $r_0$  and  $r_1$ satisfy $x_0\leq {1\over r_0}\leq {1\over r_1}\leq y_0$ or $y_1\leq {1\over r_1}\leq {1\over r_0}\leq x_1$, and if ${1\over r}= {1-\theta\over r_0}+{\theta\over r_1}$ then we have 
$$ \|f\|_{L^{p,r}}\leq C_0 \|f\|_{\dot B^{\alpha,q_0}_{r_0}}^{1-\theta}\|f\|_{\dot B^{-\beta,q_1}_{r_1}}^\theta.$$
\end{Theorem}
Finally, all these inequalities are sharp in the sense that it is not possible to remove the condition $\frac{1}{r}=\frac{1-\theta}{r_0}+\frac{\theta}{r_1}$.\\

The plan of the paper is the following: in section \ref{Definitions} we recall some facts about Lorentz and Besov spaces and we will pay a special attention to the interpolation definition of Lorentz spaces. In section \ref{Proofs} we will give the proofs and finally, in section \ref{Sharpness} we will treat the part 3) of theorem \ref{Theorem3} where we will study the sharpness of these inequalities.
\section{Functional spaces and real interpolation method}\label{Definitions}

For Besov and Triebel-Lizorkin spaces we will use the characterization based on a Littlewood-Paley decomposition. We start with a nonnegative function $\varphi\in \mathcal{D}(\mathbb{R}^n)$ such that $\varphi(\xi)=1$ over $|\xi|\leq 1/2$ and $\varphi(\xi)=0$ if $|\xi|\geq 1$. Let $\psi$ be defined as $\psi(\xi)=\varphi(\xi/2)-\varphi(\xi)$. We define the operators $S_j$ and $\Delta_j$ in the Fourier level by the formulas
$\widehat{S_j f}(\xi)=\varphi(\xi)\widehat{f}(\xi)$ and $\widehat{\Delta_j f}(\xi)=\psi(\xi)\widehat{f}(\xi)$. The distribution $\Delta_j f$ is called the $j$-th dyadic block of the Littlewood-Paley decomposition of $f$. If $\underset{j\to -\infty}{\lim} S_jf=0$ in $\mathcal{S}'(\mathbb{R}^n)$, then the equality $f=\displaystyle{\sum_{j\in \mathbb{Z}}}\Delta_jf$ is called the homogeneous Littlewood-Paley decomposition of $f$.
\begin{Definition}
For $1\leq p,q  \leq+\infty$ and $s\in \mathbb{R}$, we define the homogeneous Besov spaces as the set of distributions  $\dot{B}^{s,p}_q(\mathbb{R}^n)=\big \{f\in \mathcal{S}'(\mathbb{R}^n): \|f\|_{\dot{B}^{s,p}_q}<+\infty\big \}$, where
\begin{equation*}
\|f\|_{\dot{B}^{s,p}_q}=\left(\sum_{j\in \mathbb{Z}}2^{jsq} \|\Delta_jf\|_{L^p}^q\right)^{1/q}
\end{equation*}
with the usual modifications when $q=+\infty$. For $1<p<+\infty$,  $1\leq q \leq +\infty$ and $s\in \mathbb{R}$, we define in the same way the homogeneous Triebel-Lizorkin spaces by $\dot{F}^{s,p}_q(\mathbb{R}^n)=\big \{f\in \mathcal{S}'(\mathbb{R}^n): \|f\|_{\dot{F}^{s,p}_q}<+\infty\big \}$ with
\begin{equation*}
\|f\|_{\dot{F}^{s,p}_q}=\left\|\left(\sum_{j\in \mathbb{Z}}2^{jsq} |\Delta_jf|^q\right)^{1/q}\right\|_{L^p}
\end{equation*}
with the usual modifications when $q=+\infty$.
\end{Definition}
Note that the quantities $\|\cdot\|_{\dot{B}^{s,p}_q}$ and $\|\cdot\|_{\dot{F}^{s,p}_q}$ are only semi-norms since for $j\in \mathbb{Z}$ we have $\Delta_j P=0$ for all polynomials $P$. \\

We turn now to Lorentz spaces which are a generalization of the Lebesgue spaces. For $(X, \mu)$ a measurable space, they are usually defined in terms of the distribution and rearrangement functions $d_f(t)$ and $f^\ast(s)$ given by the formulas
$$d_f(t)=\mu(\{x: |f(x)|\geq t\})\qquad \mbox{and}  \qquad f^\ast(s)=\inf \{t: d_f(t)\leq s\},$$
where $\mu(A)$ denotes the measure of a set $A$. Then for $1\leq p< +\infty$ and $1\leq r\leq +\infty$, the Lorentz spaces $L^{p,r}(X, \mu)$ are traditionally defined in the following way
$$L^{p,r}(X, \mu)=\big\{f:X\longrightarrow\mathbb{R}: \|f\|_{L^{p,r}}<+\infty\big\}$$
where
\begin{equation}\label{Def1}
\|f\|_{L^{p,r}}=\left(\int_{0}^{+\infty} \left(s^{\frac{1}{p}}f^\ast(s)\right)^r \frac{ds}{s}\right)^{1/r},
\end{equation}
with the usual modifications when $r=+\infty$, which corresponds to the weak-$L^p$ spaces. With this characterization is not complicated to see that we have $L^{p,p}=L^p$ and that for $r_0<r_1$ we have the embedding $L^{p,r_0}\subset L^{p,r_1}$. 

However, the previous formula is not very useful since it depends on the rearrangement function $f^\ast$ and we will use a more helpful characterization which is given in the lines below.\\

We recall now some classical results from  interpolation theory concerning the real interpolation method. See \cite{BERL} for a detailed treatment. If $A_0$ and $A_1$ are two Banach spaces which are continuously embedded into a common topological vector space $V$, if $0<\theta<1$ and $1\leq r\leq +\infty$, then the real interpolation space $[A_0,A_1]_{\theta,r}$ may be defined in the following way: $f\in [A_0, A_1]_{\theta,r}$ if and only if $f\in V$ and $f$ can be written in $V$ as $f=\displaystyle{\sum_{j\in\mathbb{Z}}} f_j$, with $f_j\in A_0\cap A_1$ and $(2^{-j\theta} \|f_j\|_{A_0})_{j\in \mathbb{Z}}\in \ell^r$, $(2^{j(1-\theta)}\|f_j\|_{A_1})_{j\in \mathbb{Z}}\in \ell^r$. This space is normed with
\begin{equation}\label{DefNorm1}
\|f\|_{[A_0,A_1]_{\theta,r}} =\inf_{f=\sum f_j}  \big(\sum_{j\in\mathbb{Z}} 2^{-j\theta r} \|f_j\|^r_{A_0}\big)^{1/r}+    \big(\sum_{j\in\mathbb{Z}} 2^{j(1-\theta) r} \|f_j\|^r_{A_1}\big)^{1/r}
\end{equation}
For $f=\displaystyle{\sum_{j\in\mathbb{Z}}} f_j$ and $\rho>0$, with $\rho\neq 1$, we have the following inequality that will be very helpful in the sequel:
\begin{equation}\label{Ineq21}
 \|f\|_{[A_0,A_1]_{\theta,r}} \leq C_{\rho,\theta,r}  \big(\sum_{j\in\mathbb{Z}} \rho^{-j\theta r} \|f_j\|^r_{A_0}\big)^{(1-\theta)/r}  \big(\sum_{j\in\mathbb{Z}} \rho^{j(1-\theta) r} \|f_j\|^r_{A_1}\big)^{\theta/r}.
\end{equation}
An important property of the real interpolation method is the reiteration theorem: 

\begin{Proposition}
\begin{enumerate}
\item[]
\item[1)] If $\theta_0\neq \theta_1$, we have
\begin{equation}\label{Ineq22}
\big[ [A_0,A_1]_{\theta_0,r_0},[A_0,A_1]_{\theta_1,r_1}\big]_{\theta,r}=[A_0,A_1]_{(1-\theta)\theta_0+\theta\theta_1,r}.
\end{equation}
\item[2)] If $\theta_0=\theta_1$, (\ref{Ineq22}) is still valid if ${1\over r}={1-\theta\over r_0}+{\theta\over r_1}$.
\end{enumerate}
\end{Proposition}
We saw with the expression (\ref{Def1}) how to define Lorentz spaces $L^{p,r}(X,\mu)$ for $1<p<+\infty$, $1\leq r\leq +\infty$ as an Orlicz space. However, it will be simpler to use their characterization as real interpolates of Lebesgue spaces \cite{BERL}:
\begin{Proposition}[Lorentz spaces as interpolation spaces]\label{Proposition1}
\begin{itemize}
\item[]
\item[1)] For $1<p<+\infty$, $1\leq r\leq +\infty$
\begin{equation}\label{Ineq23}
L^{p,r}=[L^1,L^\infty]_{\theta,r}\hspace{3cm} \hbox{ with } \theta=1-{1\over p}.
\end{equation}
\item[2)] For $p_0\neq p_1$, we have
\begin{equation}\label{Ineq24}
[L^{p_0},L^{p_1}]_{\theta,r}=[L^{p_0,r_0},L^{p_1,r_1}]_{\theta,r}=L^{p,r}\qquad\hbox{ with } {1\over p}={1-\theta\over p_0}+{\theta\over p_1}.
\end{equation}
\item[3)]In the case $p_0=p_1=p$ we have
\begin{equation}\label{Ineq25}
[L^{p,r_0},L^{p,r_1}]_{\theta,r}=L^{p,r}\hspace{3cm}\hbox{ if  }\quad {1\over r}={1-\theta\over r_0}+{\theta\over r_1}.
\end{equation}
\end{itemize}
\end{Proposition}
Of course, (\ref{Ineq24}) and (\ref{Ineq25}) are consequences of (\ref{Ineq23}) through the reiteration theorem. In this paper, we shall use decompositions (\ref{DefNorm1}) and estimates (\ref{Ineq21}) when we deal with functions in Lorentz spaces, and we will mainly consider the cases $(X,\mu)=(\mathbb{R}^n,\lambda)$ where $\lambda$ is the Lebesgue measure, or $(X,\mu)=(\mathbb{Z},\mu)$ with $\mu$ the counting measure.\\

In the case of Lorentz spaces, we can use decomposition (\ref{DefNorm1}) with an useful extra property (see \cite{LEMA}) :
\begin{Lemma}\label{Lemma1}
Let  $1<p<+\infty$, $1\leq r\leq +\infty$. Then there exists a constant $C_0$ such that every $f\in L^{p,r}(X,\mu)$ can be decomposed as $f= \displaystyle{\sum_{j\in\mathbb{Z}}} f_j$ where 
\begin{itemize}
\item  $\big\|(2^{-j(p-1)/p} \|f_j\|_{L^1})\big\|_{\ell^r} + \big\|(2^{j/p}\|f_j\|_{L^\infty})\big\|_{\ell^r}\leq C_0 \|f\|_{L^{p,r}}$\\
\item the $f_j$ have disjoint supports : if $j\neq k$, $f_jf_k=0$.
\end{itemize}
\end{Lemma}
Inequality (\ref{Ineq21}) is very useful to provide an upper bound for the Lorentz norm of $f$ and it will be systematically used here. In order to get a lower bound, we shall use the following duality result : 
 
\begin{Lemma}\label{Lemma2}
Let  $1<p<+\infty$, $1\leq r\leq +\infty$. Then there exists a constant $C_0$ such that for every $f\in L^{p,r}(X,\mu)$ and every $g\in L^{\frac{p}{p-1}, \frac{r}{r-1}}(X,\mu)$, we have $fg\in L^1(X,\mu)$ and 
\begin{equation*}
\left| \int fg\ d\mu\right|\leq C_0\|f\|_{L^{p,r}} \|g\|_{L^{\frac{p}{p-1}, \frac{r}{r-1}}}.
\end{equation*}
\end{Lemma}
\begin{Remark} As we shall see in the proofs given in the section below, the characterization of Lorentz and Besov spaces based on real interpolation is a useful tool since the problem we are dealing with can be studied in terms of weighted sequences. See \cite{BERL}, \cite{DEVPOV} or \cite{KREP} for more details concerning the interpolation of Besov spaces.
\end{Remark}
\section{Refined inequalities: the proofs}\label{Proofs}
\textit{\textbf{Proof of theorem \ref{Theorem1}.}} We just write $A_\alpha(x)=\underset{j\in\mathbb{Z}}{\sup}\; 2^{j\alpha} \vert \Delta_j(x)\vert$ and $A_\beta(x)=\underset{j\in\mathbb{Z}}{\sup}\; 2^{-j\beta} \vert \Delta_j(x)\vert$ to obtain
$$ \vert f(x)\vert\leq \sum_{j\in\mathbb{Z}} \vert \Delta_jf(x)\vert\leq \sum_{j\in\mathbb{Z}} \min\bigg(2^{-j\alpha} A_{\alpha}(x),2^{j\beta}A_{\beta}(x) \bigg).$$
We define $j_0(x)$ as the largest index such that $2^{j\beta}A_{\beta}(x) \leq 2^{-j\alpha}A_{\alpha}(x) $ and we write
\begin{equation}\label{Ineq31}
\vert f(x)\vert \leq  \sum_{j\leq j_0(x)} 2^{j\beta}A_{\beta}(x)  +\sum_{j>j_0(x)} 2^{-j\alpha}A_{\alpha}(x) \leq  C A_{\alpha}(x)^{\beta\over \alpha+\beta}  A_{\beta}(x)^{\alpha\over \alpha+\beta},
\end{equation}
thus, inequality (\ref{Ineq5}) is proved. In order to obtain (\ref{Ineq6}), it is enough to apply H\"older inequality in the expression above since we have $\theta=\frac{\alpha}{\alpha+\beta}$ and $\frac{1}{p}=\frac{1-\theta}{q_0}+\frac{\theta}{q_1}$.\hfill$\blacksquare$
\begin{Remark}
Inequality (\ref{Ineq31}) is a little more precise than Hedberg's inequality \cite{ADAH, HEDB} : if $f=I_{\alpha}g(x)$ where $I_{\alpha}$ is a Riesz potential\footnote{defined in the Fourier level by $\widehat{I_{\alpha}g}(\xi)= \vert\xi\vert^{-\alpha} \widehat{g}(\xi)$.} with $0<\alpha<n$ and if $g\in \dot B^{-\beta,\infty}_\infty$, then, if $\mathcal{M}_g$ is the Hardy--Littlewood maximal function of $g$, we have $A_{\alpha}(x)\leq C \mathcal{M}_g(x)$ and $A_{\beta}(x)\leq C \|g\|_{ \dot B^{-\beta,\infty}_{\infty}}$. Thus, we find easily the refined Sobolev inequality (\ref{Ineq1}). See more details in \cite{ ORU, CHAM}.\\
\end{Remark}
\textit{\textbf{Proof of theorem \ref{Theorem2}.}} We start picking $p_0$ and $p_1$ such that $1\leq q_0<p_0<p<p_1<q_1\leq +\infty$ with $\frac{2}{p}= \frac{1}{p_0}+\frac{1}{p_1}$. We have then $\frac{1}{p_i}=\frac{1-a_i}{q_0}+\frac{a_i}{q_1}$ with $0<a_i<1$ and $i=0,1$. We write
$$\|\Delta_jf\|_{L^{p_i}}\leq\|\Delta_jf\|_{L^{q_0}}^{1-a_i}\|\Delta_jf\|_{L^{q_1}}^{a_i}= \left(2^{j\alpha}\|\Delta_jf\|_{L^{q_0}}\right)^{1-a_i}\left( 2^{-j\beta} \|\Delta_jf\|_{L^{q_1}}\right)^{a_i} 2^{j\big[-\alpha (1-a_i)+\beta a_i\big]}.$$
Recalling that $\frac{1}{p}=\frac{1-\theta}{q_0}+\frac{\theta}{q_1}$ and $\theta=\frac{\alpha}{\alpha+\beta}$ we have $-\alpha(1-a_0)+\beta a_0=\alpha(1-a_1)-\beta a_1$. Thus, noting $\rho=2^{-2[\alpha(1-a_0)-\beta a_0]}>0$ and using the H\"older inequality we obtain 
\begin{eqnarray*}
\sum_{j\in \mathbb{Z}}\rho^{-jr/2}\|\Delta_jf\|^r_{L^{p_0}}\leq \|f\|_{\dot{B}^{\alpha, q_0}_r}^{r(1-a_0)}\|f\|_{\dot{B}^{-\beta, q_1}_r}^{ra_0 } &\mbox{ and \quad} & \sum_{j\in \mathbb{Z}}\rho^{jr/2}\|\Delta_jf\|^r_{L^{p_1}}\leq  \|f\|_{\dot{B}^{\alpha, q_0}_r}^{r(1-a_1)}\|f\|_{\dot{B}^{-\beta, q_1}_r}^{ra_1}.\\
\end{eqnarray*}
From this, and applying proposition \ref{Proposition1}, we deduce that if $f\in \dot{B}^{\alpha, q_0}_r\cap \dot{B}^{-\beta, q_1}_r(\mathbb{R}^n)$ then $f\in [L^{p_0}, L^{p_1}]_{\frac{1}{2},r}=L^{p,r}$. Furthermore, using inequality (\ref{Ineq21}) we finally have:
$$\|f\|_{L^{p,r}}\leq C_{\rho,r} \|f\|_{\dot{B}^{\alpha, q_0}_r}^{1-\theta}\|f\|_{\dot{B}^{-\beta, q_1}_r}^{\theta}$$
\hfill$\blacksquare$\\ 
\textit{\textbf{Proof of theorem \ref{Theorem3}.}}
\paragraph{1) Case $r>r^*$: } With no loss of generality, we may assume that $q_0<q_1$ and we fix $\varepsilon>0$ such that 
$${1\over q_1}<{1\over p}-\varepsilon({1\over q_0}-{1\over q_1})={1\over p_1}<{1\over p}+\varepsilon({1\over q_0}-{1\over q_1})={1\over p_0}<{1\over q_0}.$$
The proof follows essentially the same ideas used in the previous theorem. Indeed, we have, for $\gamma_j=2^{j\alpha} \|\Delta_jf\|_{L^{q_0}}$ and $\eta_j=2^{-j\beta} \|\Delta_jf\|_{L^{q_1}}$, and for $\epsilon_0=1$ and $\epsilon_1=-1$,  
$$ \|\Delta_j f\|_{L^{p_i}}\leq \|\Delta_jf\|_{L^{q_0}}^{1-\theta+\epsilon_i\varepsilon} \|\Delta_jf\|_{L^{q_1}}^{\theta-\epsilon_i\varepsilon}= \gamma_j^{1-\theta+\epsilon_i\varepsilon} \eta_j^{\theta-\epsilon_i\varepsilon}2^{-j\epsilon_i\varepsilon(\alpha+\beta)}.$$
As $r_0\neq r_1$, we can only say that 
$(\gamma_j^{1-\theta+\epsilon_i\varepsilon} \eta_j^{\theta-\epsilon_i\varepsilon})_{j\in\mathbb{Z}}\in \ell^{\rho_i}$ where ${1\over\rho_i}={1-\theta+\epsilon_i\varepsilon\over r_0}+{\theta-\epsilon_i\varepsilon\over r_1}$. 
We may use inequality (\ref{Ineq21}), but we get only that $f\in [L^{p_0},L^{p_1}]_{1/2,\rho}=L^{p,\rho}$  with $\rho=\max(\rho_0,\rho_1)$,   and satisfies inequality (\ref{Ineq9}) with $r=\rho$. However, we may choose $\varepsilon$ as small as we want, and thus $\rho$ as close to $r^*$ as we want; thus $f$ satisfies (\ref{Ineq9}) for every $r>r^*$.
\paragraph{2) Case $r=r^*$: }
\begin{itemize}
\item[a)] if $r=r_0=r_1$ : this case was treated in theorem \ref{Theorem2}.
\item[b)] if $r_0=q_0$ and $r_1=q_1$ : This is a  direct consequence of (\ref{Ineq6}) since we have $\|f\|_{\dot B^{\alpha,q_i}_{q_i}}=\|f\|_{\dot F^{\alpha,q_i}_{q_i}}\subset  \|f\|_{\dot F^{\alpha,q_i}_{\infty}}$ and $\|f\|_{\dot B^{-\beta,q_i}_{q_i}}=\|f\|_{\dot F^{-\beta,q_i}_{q_i}}\subset  \|f\|_{\dot F^{-\beta,q_i}_{\infty}}$, we obtain (\ref{Ineq7}).
\item[c)]  Case $1<p\leq 2$ and $r^*=p$: We just write
$$ \|\Delta_jf\|_{L^p}\leq \|\Delta_jf\|_{L^{q_0}}^{1-\theta} \|\Delta_jf\|_{L^{q_1}}^{\theta} =(2^{\alpha} \|\Delta_jf\|_{L^{q_0}})^{1-\theta} (2^{-j\beta} \|\Delta_jf\|_{L^{q_1}})^{\theta}$$
and get by H\"older inequality:
$$ \|f\|_{\dot B^{0,p}_p}\leq C \|f\|_{\dot B^{\alpha,q_0}_{r_0}}^{1-\theta}\|f\|_{\dot B^{-\beta,q_1}_{r_1}}^\theta.$$
We then use  the embedding $\dot B^{0,p}_{p}\subset L^p=L^{p,p}$, which is valid for $p\leq 2$.\hfill$\blacksquare$\\ 
\end{itemize}
\textit{\textbf{Proof of theorem \ref{Theorem4}.}} We see how direct interpolation has given us theorem \ref{Theorem3}, but we only obtained partial results for theorems \ref{Theorem4} and \ref{Theorem5}. Indeed, in theorem \ref{Theorem4}, we want a positive result for $q_0\leq r_0\leq r_1\leq q_1$ but thus far we have proven the result only for $(r_0,r_1)=(q_0,q_1)$ and for $q_0\leq r_0=r_1\leq q_1$. To complete the proof of theorem \ref{Theorem4}, we must reiterate interpolations to those new estimates.\\

This will be done through the following lemma :

\begin{Lemma}\label{Lemma3}
Let $\alpha, \beta>0$, $q_0,q_1\in [1+\infty]$, $q_0< q_1$. Let $\theta={\alpha\over \alpha+\beta}\in ]0,1[$ and let ${1\over p}={1-\theta\over q_0}+{\theta\over q_1}$.
\begin{itemize}
\item[1)]  If $q_0\leq r_0\leq q_1$ and let ${1\over r}= {1-\theta\over r_0}+{\theta\over q_1}$, then we have 
\begin{equation}\label{IneqLemma1}
\|f\|_{L^{p,r}}\leq C_0 \|f\|_{\dot B^{\alpha,q_0}_{r_0}}^{1-\theta}\|f\|_{\dot B^{-\beta,q_1}_{q_1}}^\theta.
\end{equation}
\item[2)]  If $q_0\leq r_1\leq q_1$ and let ${1\over r}= {1-\theta\over q_0}+{\theta\over r_1}$, then:
\begin{equation}\label{IneqLemma2}
\|f\|_{L^{p,r}}\leq C_0 \|f\|_{\dot B^{\alpha,q_0}_{q_0}}^{1-\theta}\|f\|_{\dot B^{-\beta,q_1}_{r_1}}^\theta.
\end{equation}
\end{itemize}
\end{Lemma}
\textit{\textbf{Proof of the lemma \ref{Lemma3}.}} We only prove the first inequality, as the proof for the second one is similar. Since $f\in \dot{B}^{\alpha,q_0}_{r_0}$, noting $\lambda_j=2^{j\alpha}\|\Delta_jf\|_{L^{q_0}} $ we have  $(\lambda_j)_{j\in\mathbb{Z}}\in \ell^{r_0}$.  Thus, using lemma \ref{Lemma1} for the interpolation 
\begin{equation}\label{Interpolation}
\ell^{r_0}=[ \ell^{q_0},\ell^{q_1}]_{\eta,r},
\end{equation}
with $\frac{1}{r_0}=\frac{1-\eta}{q_0}+\frac{\eta}{q_1}$, we see that we have a partition $\mathbb{Z}=\displaystyle{\sum_{k\in\mathbb{Z}}} Z_k$ such that
\begin{equation}\label{Ineq32}
\bigg\|2^{-k\eta}  \big(\sum_{j\in Z_k}\!  \lambda_j^{q_0}\big)^{1\over q_0}\bigg\|_{\ell^{r_0}}  + \bigg\|2^{k(1-\eta)}  \big(\sum_{j\in Z_k}\!  \lambda_j^{q_1}\big)^{1\over q_1}\bigg\|_{\ell^{r_0}} \! \leq \!  C \|\lambda_j\|_{\ell^{r_0}}
\end{equation}
Moreover since $f\in \dot B^{-\beta,q_1}_{q_1}$ we have 
$$ \bigg( \big(\sum_{j\in Z_k} 2^{-j\beta q_1} \|\Delta_jf \|_{L^{q_1}}^{q_1}\big)^{1/q_1} \bigg)_{k\in\mathbb{Z}}\in \ell^{q_1}.$$
Let us note $\alpha_k= \big(\displaystyle{\sum_{j\in Z_k}}2^{-j\beta q_1} \|\Delta_jf \|_{L^{q_1}}^{q_1}\big)^{1/q_1}$, $\beta_k=  2^{-k\eta}  \big( \displaystyle{\sum_{j\in Z_k}}\lambda_j^{q_0}\big)^{1\over q_0}$, $\gamma_k= 2^{k(1-\eta)}  \big(\displaystyle{\sum_{j\in Z_k}}  \lambda_j^{q_1})^{1\over q_1}$ and $f_k=\displaystyle{\sum_{j\in Z_k}} \Delta_j f$. 
We apply now inequality (\ref{Ineq7}) and theorem \ref{Theorem2} to obtain
\begin{equation}\label{Ineq33}
\|f_k\|_{L^{p}}\leq  C_A \|f_k\|_{\dot B^{\alpha,q_0}_{q_0}}^{1-\theta}\|f_k\|_{\dot B^{-\beta,q_1}_{q_1}}^\theta\leq C  \alpha_k^{\theta}\beta_k^{1-{\theta}} 2^{k\eta(1-{\theta})}
\end{equation}
and
\begin{equation}\label{Ineq34}
\|f_k\|_{L^{p,q_1}}\leq  C_B \|f_k\|_{\dot B^{\alpha,q_0}_{q_1}}^{1-\theta}\|f_k\|_{\dot B^{-\beta,q_1}_{q_1}}^\theta\leq C  \alpha_k^{\theta}\gamma_k^{1-{\theta}} 2^{-k(1-\eta)(1-{\theta})}.
\end{equation}
Since we have $f=\displaystyle{\sum_{k\in\mathbb{Z}}} f_k$, with these two inequalities at hand, and using (\ref{Ineq21}), we find that $f\in  [L^{p},L^{p,q_1}]_{\eta,r}$ with  $\frac{1}{r}=\frac{1-\eta}{p}+\frac{\eta}{q_1}$. But, since $\frac{1}{r_0}=\frac{1- \eta}{q_0}+\frac{\eta}{q_1}$ and $\frac{1}{p}=\frac{1-\theta}{q_0}+\frac{\theta}{q_1}$, we obtain $[L^{p},L^{p,q_1}]_{\eta,r}=L^r$ with $\frac{1}{r}=\frac{1-\theta}{r_0}+\frac{\theta}{q_1}$. \hfill$\blacksquare$\\

\begin{Remark}
Note that we use twice interpolation arguments: first in estimate (\ref{Ineq32}) and then with inequalities (\ref{Ineq33}) and (\ref{Ineq34}) in order to obtain $f\in  [L^{p},L^{p,q_1}]_{\eta,r}$.\\
\end{Remark}
Once this lemma is proved, it is enough to reapply similar arguments to obtain theorem \ref{Theorem4}. Indeed, since we have $q_0<r_0<r_1<q_1$, we start using $\ell^{r_0}=[ \ell^{q_0}, \ell^{r_1}]_{\eta, {r_0}}$ instead of (\ref{Interpolation}) and we obtain a partition $\mathbb{Z}=\displaystyle{\sum_{k\in\mathbb{Z}}} Z_k$ such that
\begin{equation*}
\bigg\|2^{-k\eta}  \big(\sum_{j\in Z_k}\!  \lambda_j^{q_0}\big)^{1\over q_0}\bigg\|_{\ell^{r_0}}  + \bigg\|2^{k(1-\eta)}  \big(\sum_{j\in Z_k}\!  \lambda_j^{r_1}\big)^{1\over r_1}\bigg\|_{\ell^{r_0}} \! \leq \!  C \|\lambda_j\|_{\ell^{r_0}}
\end{equation*}
with $\frac{1}{r_0}=\frac{1-\eta}{q_0}+\frac{\eta}{r_1}$ and where the sequence $(\lambda_j)_{j\in \mathbb{N}}$ with $\lambda_j=2^{j\alpha}\|\Delta_jf\|_{L^{q_0}}$ belongs to $\ell^{r_0}$ since $f\in \dot{B}^{\alpha, q_0}_{r_0}$.\\
Since $f\in \dot B^{-\beta,q_1}_{r_1}$ we have $ \big(\big(\sum_{j\in Z_k} 2^{-j\beta q_1} \|\Delta_jf \|_{L^{q_1}}^{q_1}\big)^{1/q_1} \big)_{k\in\mathbb{Z}}\in \ell^{q_1}$, and we note again $\alpha_k= \big(\displaystyle{\sum_{j\in Z_k}}2^{-j\beta q_1} \|\Delta_jf \|_{L^{q_1}}^{q_1}\big)^{1/q_1}$, $\beta_k=  2^{-k\eta}  \big( \displaystyle{\sum_{j\in Z_k}}\lambda_j^{q_0}\big)^{1\over q_0}$, $\gamma_k= 2^{k(1-\eta)}  \big(\displaystyle{\sum_{j\in Z_k}}  \lambda_j^{q_1})^{1\over q_1}$ and $f_k=\displaystyle{\sum_{j\in Z_k}} \Delta_j f$. 
Next, we only need to apply (\ref{IneqLemma2}) and (\ref{Ineq8}) instead of  (\ref{Ineq33}) and (\ref{Ineq34}) to obtain
\begin{equation*}
\|f_k\|_{L^{p,s}}\leq  C_A \|f_k\|_{\dot B^{\alpha,q_0}_{q_0}}^{1-\theta}\|f_k\|_{\dot B^{-\beta,q_1}_{r_1}}^\theta\leq C  \alpha_k^{\theta}\beta_k^{1-{\theta}} 2^{k\eta(1-{\theta})}
\end{equation*}
where $\frac{1}{s}=\frac{1-\theta}{q_0}+\frac{\theta}{r_1}$, and
\begin{equation*}
\|f_k\|_{L^{p,r_1}}\leq  C_B \|f_k\|_{\dot B^{\alpha,q_0}_{r_1}}^{1-\theta}\|f_k\|_{\dot B^{-\beta,q_1}_{r_1}}^\theta\leq C  \alpha_k^{\theta}\gamma_k^{1-{\theta}} 2^{-k(1-\eta)(1-{\theta})}.
\end{equation*}
Finally, we have via inequality (\ref{Ineq21}) that $f\in  [L^{p,s},L^{p,r_1}]_{\eta,r}$ with  $\frac{1}{r}=\frac{1-\eta}{s}+\frac{\eta}{r_1}$. To conclude, we use the fact that $\frac{1}{s}=\frac{1-\theta}{q_0}+\frac{\theta}{r_1}$ and $\frac{1}{r_0}=\frac{1-\eta}{q_0}+\frac{\eta}{r_1}$ in order to obtain that $f\in L^{p,r}$ with $\frac{1}{r}=\frac{1-\theta}{r_0}+\frac{\theta}{r_1}$ . \hfill$\blacksquare$\\[5mm]
\textit{\textbf{Proof of theorem \ref{Theorem5}.}}  Similarly, in theorem \ref{Theorem5}, we want a positive result on the triangles $1/y_0\leq r_1\leq r_0\leq 1/x_0$ and $1/x_1\leq r_0\leq r_1\leq 1/y_1$,  and we have already obtained that the theorem is true for $1/y_0\leq r_1= r_0\leq 1/x_0$ and $1/x_1\leq r_0= r_1\leq 1/y_1$ as well as for $(1/x_0,1/y_0)$ and $(1/x_1,1/y_1)$.  To complete the proof of theorem \ref{Theorem5}, we must reiterate interpolations to those new estimates. This will be achieved with the following lemma. 

\begin{Lemma}\label{Lemma4}
Let $\Sigma$ be the set of points $(x,y)\in [0,1]\times [0,1]$ such that, for $z^*=(1-\theta)x+\theta y$, we have the inequality
$$\|f\|_{L^{p,1/z^*}}\leq C_{x,y}\|f\|_{\dot B^{\alpha_0,q_0}_{1/x}}^{1-\theta}\|f\|_{\dot B^{-\beta,q_1}_{1/y}}^\theta.$$ 
and let $A,B,C\in [0,1]\times [0,1]$ such that $[A,B]$ is horizontal $(y_A=y_B)$ and  $[A,C]$ is vertical ($x_A=x_C$). Then :
\begin{itemize}
\item if $A\in \Sigma$ and $B\in \Sigma$ then $[A,B]\subset \Sigma$,
\item if $A\in \Sigma$ and $C\in \Sigma$ then $[A,C]\subset \Sigma$,
\item if $A\in \Sigma$ and $[B,C]\subset \Sigma$, then the triangle $ABC$ is contained in $\Sigma$.
\end{itemize}
\end{Lemma}
\textit{\textbf{Proof of the lemma \ref{Lemma4}.}} The proof of this inequality follows closely the ideas of lemma \ref{Lemma3}. We give the details here for the sake of completness.\\ 

We begin with the case of $M\in[A,B]$: set $x_A, z_0, \rho$ such that $z_0=(1-\theta)x_A+\theta \rho$ and $x_B, z_1, \rho$ such that $z_1=(1-\theta)x_B+\theta \rho$. Then, if we define $x_M$ such that $x_M=(1-\eta)x_A+\eta x_B$, we must show that we have the inequality
$$\|f\|_{L^{p,1/z}}\leq C_0 \|f\|_{\dot B^{\alpha,q_0}_{1/x_M}}^{1-\theta}\|f\|_{\dot B^{-\beta,q_1}_{1/\rho}}^\theta$$
with $z=(1-\theta)x_M+\theta \rho$. 

\begin{itemize}
\item Since $f\in \dot B^{\alpha,q_0}_{1/x_M}$, we have that $(\lambda_j)_{j\in \mathbb{Z}}\in \ell^{1/x_M}$, with $\lambda_j=2^{j\alpha}\|\Delta_j f\|_{L^{q_0}}$. Moreover, by hypothesis we have $x_M=(1-\eta)x_A+\eta x_B$, so we can write $\ell^{1/x_M}=[\ell^{1/x_A}, \ell^{1/x_B}]_{\eta, 1/x_M}$. Thus, with lemma \ref{Lemma1} we obtain a partition $\mathbb{Z}=\displaystyle{\sum_{k\in\mathbb{Z}}} Z_k$ such that 
\begin{equation*}
\bigg\|2^{-k\eta}  \big(\sum_{j\in Z_k}\lambda_j^{1/x_A}\big)^{x_A}\bigg\|_{\ell^{1/x_M}}+\bigg\|2^{k(1-\eta)}\big(\sum_{j\in Z_k} \lambda_j^{1/x_B}\big)^{x_B}\bigg\|_{\ell^{1/X_M}} \leq C \|\lambda_j\|_{\ell^{1/x_M}}.
\end{equation*}
We will note $\beta_k=2^{-k\eta}  \big(\displaystyle{\sum_{j\in Z_k}}\lambda_j^{1/x_A}\big)^{x_A}$,  $\gamma_k=2^{k(1-\eta)}\big(\displaystyle{\sum_{j\in Z_k}} \lambda_j^{1/x_B}\big)^{x_B}$ and $f_k=\displaystyle{\sum_{j\in Z_k}}\Delta_jf$. 
\item Since  $f\in \dot B^{-\beta,q_1}_{1/\rho}$, we can write $\alpha_k=\bigg( \big(\displaystyle{\sum_{j\in Z_k}} 2^{-j\beta 1/\rho} \|\Delta_jf \|_{L^{q_1}}^{1/\rho}\big)^{\rho} \bigg)_{k\in\mathbb{Z}}\in \ell^{1/\rho}.$
\end{itemize}
Now, since $z_0=(1-\theta)x_A+\theta \rho$ and $z_1=(1-\theta)x_B+\theta \rho$ we can apply theorem \ref{Theorem4} to the functions $f_k$ to obtain:
\begin{eqnarray*}
\|f_k\|_{L^{p, 1/z_0}} &\leq  &C_A \|f_k\|_{\dot B^{\alpha,q_0}_{1/x_A}}^{1-\theta}\|f_k\|_{\dot B^{-\beta,q_1}_{1/\rho}}^\theta\leq C  \alpha_k^{\theta}\beta_k^{1-{\theta}} 2^{k\eta(1-{\theta})}\\[5mm]
\|f_k\|_{L^{p, 1/z_1}} &\leq  &C_B \|f_k\|_{\dot B^{\alpha,q_0}_{1/x_B}}^{1-\theta}\|f_k\|_{\dot B^{-\beta,q_1}_{1/\rho}}^\theta\leq C  \alpha_k^{\theta}\gamma_k^{1-{\theta}} 2^{-k(1-\eta)(1-{\theta})}
\end{eqnarray*}
From these two estimates we deduce, since $f=\displaystyle{\sum_{k\in \mathbb{Z}}}f_k$, that $f\in [L^{p,1/z_0}, L^{p,1/z_1}]_{\eta, 1/z}$ with $z=(1-\eta)z_0+\eta z_1$. But by hypothesis we have $z_0=(1-\theta)x_A+\theta \rho$ and $z_1=(1-\theta)x_B+\theta \rho$, so we find that $z=(1-\theta)x_M +\theta \rho$. \\

We have proven that we may interpolate along horizontal lines and the vertical case is totally similar. To finish the proof of the lemma suppose that $A\in \Sigma$ and $[B,C]\subset \Sigma$ and take a point $P\in \Sigma$. Write $P\in [M,N]$ where $[M,N]$ is horizontal, $M \in [A,C]$ and $N\in [B,C]$.  Then we find that $M\in\Sigma$ by vertical interpolation between $A$ and $C$, and that $P\in\Sigma$ by horizontal interpolation between $M$ and $N$.\\

Thus, Lemma \ref{Lemma4} is proven, and this finishes the proof of theorem \ref{Theorem5}.  \hfill$\blacksquare$\\

\section{Sharpness and optimality of the inequalities.}\label{Sharpness}
In this section, we adapt Bahouri and Cohen's example \cite{COHB} to the general case. Their idea is to use the analysis of the regularity of chirps by Jaffard and Meyer \cite{JAFM}. They express their example in terms of wavelets, as it is easy to estimate Besov norms in wavelet bases \cite{LEMM, MEY}. They pick just one wavelet per scale and adjust its support so that they are able to compute the Lorentz norm of the sum.  As we shall see, their example can be extended to the general case but we will not use wavelets, but atoms for Besov spaces \cite{FRAJ, TRI}, as we shall use only upper estimates for Besov norms. We consider a function $\omega$ such that it is supported in the ball $B(0,1)$, is ${\cal C}^N(\mathbb{R}^n)$  and $\displaystyle{\int_{\mathbb{R}^n}} x^\gamma\omega(x) \ dx=0$ for all $\gamma\in\mathbb{N}^n$  with $\vert\gamma\vert<N$, for some $N$ such that $N> \max(\vert \alpha\vert, \vert \beta\vert )$. Then we have, for all $1\leq q\leq+\infty$, $1\leq r\leq +\infty$ and $\vert s\vert<N$,
\begin{equation}\label{Sharp1}
\bigg\| \sum_{j\in\mathbb{Z}}\sum_{k\in\mathbb{Z}^n} \lambda_{j,k}\omega(2^jx-k)\bigg\|_{\dot B^{s,q}_r}\leq C \sum_{j\in\mathbb{Z}} 2^{j s r} 2^{-jnr/q} \big(\! \sum_{k\in\mathbb{Z}^n} \vert\lambda_{j,k}\vert^q\big)^{r/q}\big)^{1/r}.
\end{equation}
For our example, we shall fix some $X\in\mathbb{R}$, some $Y\in\mathbb{R}$ and some $\gamma\in \mathbb{R}$ and we define 
$$f_L=\sum_{j=1}^{j_L} 2^{jX} \sum_{ k\in K_j} \omega(2^jx-k)\qquad \mbox{and}\qquad g_L= \sum_{j=j_1}^{j_L} 2^{jY} \sum_{ k\in K_j} \omega(2^jx-k)$$
where the $K_j$ are finite sets such that
\begin{itemize}
\item the supports of the functions $(\omega(2^jx-k))_{j_1\leq j\leq j_L,k\in K_j}$ are disjoint each one from each other
\item $K_j$ is a set of cardinal $j$ with $ 2^{\delta j}\leq A_j\leq \ 2^{\delta (j+1)}$ (which is possible if $\delta j_1\geq 0$  and $\delta j_L\geq 0$).
\end{itemize}
For every $r\in [1,+\infty]$, we then have the following estimates (from (\ref{Sharp1})):
\begin{eqnarray*}
\|f_L\|_{\dot B^{\alpha,q_0}_{r}}&\leq &C_r  \Big(\sum_{j=j_1}^{j_L}  2^{jr (X+\alpha -n/q_0+\delta/q_0)} \Big)^{1/r}\\
\|f_L\|_{\dot B^{-\beta,q_1}_{r}}&\leq &C_r  \Big(\sum_{j=j_1}^{j_L}  2^{jr (X- \beta -n/q_1+\delta/q_1)} \Big)^{1/r}\\
\|g_L\|_{\dot B^{-\alpha,q_0/(q_0-1)}_r}&\leq &C_r  \Big(\sum_{j=j_1}^{j_L}  2^{jr (Y-\alpha -n(1-1/q_0)+\delta(1-1/q_0))} \Big)^{1/r}\\
\|g_L\|_{\dot B^{\beta,q_1/(q_1-1)}_r}&\leq &C_r   \Big(\sum_{j=j_1}^{j_L}  2^{jr (Y+\beta -n(1-1/q_1)+\delta(1-1/q_1))} \Big)^{1/r}
\end{eqnarray*}
We thus fix $X,Y,\gamma$ such that
\begin{equation}\label{Cas1}
\begin{cases} X+\alpha -n/q_0+\delta/q_0=0\\
X-\beta-n/q_1+\delta/q_1=0\\
Y-\alpha -n(1-1/q_0)+\delta(1-1/q_0)=0
\end{cases}
\end{equation}
Since $\alpha-n/q_0+\delta/q_0= -\beta-n/q_1+\delta/q_1$, we have as well $Y+\beta-n(1-1/q_1)+\delta(1-1/q_1)=0$. This gives that, for all $r\in [1,+\infty]$ we have
\begin{eqnarray*}
\|f_L\|_{\dot B^{\alpha,q_0}_r}&\leq & C_r L^{1/r} \\ 
\|f_L\|_{\dot B^{-\beta,q_1}_r}&\leq & C_r L^{1/r} \\
\|g_L\|_{\dot B^{-\alpha,q_0/(q_0-1)}_r}&\leq & C_r L^{1/r}\\
\|g_L\|_{\dot B^{\beta, q_1/(q_1-1)}_r}&\leq & C_r L^{1/r}
\end{eqnarray*}
where $C_r$ does not depend on $L$. Theorem \ref{Theorem2} shows that
$$\| f_L\|_{L^{p,r}}\leq C_r L^{1/r}\qquad \mbox{and}\qquad \|g_L\|_{L^{p/(p-1),r}}\leq C_r L^{1/r}.$$
Moreover, since the supports of the functions $\omega(2^jx-k)$ are disjoint, we have
$$ \int f_L\, g_L\, dx=\sum_{j=j_1}^{j_L} 2^{j(X+Y)} \sum_{ k\in K_j} 2^{-jn} \|\omega\|_{L^2}^2\geq \sum_{j=j_1}^{j_L} 2^{j(X+Y-n+\delta)}  \|\omega\|_{L^2}^2 $$
From (\ref{Cas1}), we have $X+Y-n+\delta=0$, so that, using lemma \ref{Lemma2} we obtain:
$$ \|\omega\|_{L^2}^2\  L\leq C \|f_L\|_{L^{p,r}} \|g_L\|_{L^{{p\over p-1},{r\over r-1}}}\leq C' \|f_L\|_{L^{p,r}} L^{1-{1\over r}}$$
If we asssume that the interpolation inequality (\ref{Ineq8}) is valid, we get that
 $$ L^{1/r}\leq C \|f_L\|_{L^{p,r}}\leq C'   \|f_L\|_{\dot B^{\alpha,q_0}_{r_0}}^{1-\theta}\|f_L\|_{\dot B^{-\beta,q_1}_{r_1}}^\theta\leq C'' L^{(1-\theta)/r_0+\theta/r_1}$$
Letting $L\longrightarrow +\infty$, we find that $1/r\leq (1-\theta)/r_0+\theta/r_1$ : we have thus proven the optimality of these inequalities.


\end{document}